# PERIODIC BOXCAR DECONVOLUTION AND DIOPHANTINE APPROXIMATION

By Iain M. Johnstone[1] and Marc Raimondo[2]

*Stanford University and University of Sydney*

We consider the nonparametric estimation of a periodic function that is observed in additive Gaussian white noise after convolution with a "boxcar," the indicator function of an interval. This is an idealized model for the problem of recovery of noisy signals and images observed with "motion blur." If the length of the boxcar is rational, then certain frequencies are irretrievably lost in the periodic model. We consider the rate of convergence of estimators when the length of the boxcar is *irrational*, using classical results on approximation of irrationals by continued fractions. A basic question of interest is whether the minimax rate of convergence is slower than for nonperiodic problems with $1/f$-like convolution filters. The answer turns out to depend on the type and smoothness of functions being estimated in a manner not seen with "homogeneous" filters.

**1. Introduction.**

1.1. *Statement of problem and motivation.* Suppose that we observe $Y(t)$ for $t \in [-1, 1]$, where $Y$ is drawn from an indirect estimation model in Gaussian white noise:

(1) $$Y(t) = \int_{-1}^{t} K_a f(s) \, ds + \epsilon W(t),$$

where

(2) $$K_a f(t) = \frac{1}{2a} \int_{-a}^{a} f(t - u) \, du, \qquad a > 0,$$

Received May 2002; revised July 2003.

[1]Supported in part by NSF Grants DMS-95-05151, DMS-00-72661 and NIH Grant RO1 CA 72028.

[2]Supported in part by the Australian Academy of Science and by NSF Grant DMS 96-31278.

$\{W(t), t \in [-1, 1]\}$ is a standard two-sided Wiener process and $\epsilon$ is small and assumed known. It is desired to estimate the unknown signal $f$, assumed to be periodic on $[-1, 1]$. We refer to this as boxcar deconvolution, because $K_a f = f \star k_a$ corresponds to convolution with the step function $k_a(t) = (2a)^{-1} I\{|t| \leq a\}$.

The problem has the peculiar feature that if the boxcar half-width $a$ is *rational*, then certain frequencies are completely unrecoverable from the data. Indeed, because of the periodic and convolution structure, the problem is diagonalized in the Fourier basis. Thus, let $e_k(t) = e^{\pi i k t}$, for integer $k \in \mathbb{Z}$. Then $K_a e_k = r_k e_k$, where the eigenvalues $r_0 = 1$ and

$$(3) \qquad r_k = \frac{\sin \pi k a}{\pi k a}, \qquad k \neq 0.$$

Furthermore, setting $y_k = \int_{-1}^{1} e_k(t)\, dY(t)$, $\theta_k = \langle f, e_k \rangle := \int_{-1}^{1} f(t) e_k(t)\, dt$, and $z_k = \int_{-1}^{1} e_k(t)\, dW(t)$, we find that model (1) is equivalent to

$$(4) \qquad y_k = r_k \theta_k + \epsilon z_k, \qquad k \in \mathbb{Z}.$$

For rational $a = p/q$, the eigenvalues $r_k$ vanish for all integer multiples $k = jq$ of $q$. In the Fourier expansion $\sum \langle f, e_k \rangle e_k$, all information about the coefficients $\langle f, e_{jq} \rangle$ is lost after convolution. For irrational $a$, however, the inversion formula

$$(5) \qquad \langle f, e_k \rangle = \frac{1}{r_k} \langle K_a f, e_k \rangle$$

is at least well defined, since $r_k \neq 0$ for any $k \in \mathbb{Z}$. The object of this paper is to study the quality of estimation of $f$ attainable for irrational $a$ in the small noise limit $\epsilon \to 0$.

Motivation for studying this special problem arises from several sources:

(i) It may be viewed as an idealization of the problem of recovery from linear motion blur plus noise in a fixed field of view. If a camera is passing over a scene $f(x, y)$ along a direction $(1, r)$ at unit speed, then in exposure time $2a$ the image acquired at point $(x, y)$ may be modeled as

$$(6) \qquad Kf(x, y) = \frac{1}{2a} \int_{-a}^{a} f(x + u, y + ru)\, du.$$

Our model is a one-dimensional version of horizontal motion, $r = 0$. While the periodicity assumption on $f$ may seem artificial, it does capture the property that if $f$ is locally periodic with period $2a$ near $(x, y)$ (as in certain textures), then $Kf$ is locally constant near $(x, y)$. Compare the discussion in Section 5.1. A more detailed discussion of linear motion blur, with photographic examples, may be found in Bertero and Boccacci [(1998), pages 54–58].

(ii) It is related to the problem of periodic density estimation with uniform errors. Suppose $X_1, \ldots, X_n$ are i.i.d. random variables with unknown periodic density $f$ on the circle $\mathbb{T}$. However, the $X_i$ are not observed; instead we see jittered versions

$$Y_i = X_i + z_i,$$

where $\{z_i\}$ are i.i.d. uniformly distributed on $[-a, a]$ and circular addition is used.

(iii) As an inverse problem, (5) is nonstandard: the eigenvalues $r_k$ oscillate inside an envelope decaying like $1/frequency$, for $k \neq 0$,

$$r_k \leq c/|k|, \qquad c = (\pi a)^{-1}.$$

We may ask the following: is the quality of estimation—measured by minimax rate of convergence as $\epsilon \to 0$—determined by the $1/|k|$ decay, or is it affected by the oscillatory behavior?

(iv) Let $\|x\|$ denote the distance from $x \in \mathbb{R}$ to the nearest integer. For $k \neq 0$,

$$(7) \qquad \frac{2}{\pi} \frac{\|ka\|}{|ka|} \leq r_k \leq \frac{\|ka\|}{|ka|},$$

and so the oscillations in (3) are driven by

$$(8) \qquad \|ka\| := \inf\{|ka - l|, l \in \mathbb{Z}\}.$$

The study of such "Diophantine approximations" uses the classical theory of continued fractions, for example, Lang (1966) and Khinchin (1992), and plays a basic role in this paper.

There is a large literature on statistical inverse problems—for some recent reviews see Tenorio (2001) and Evans and Stark (2002). In particular, the sequence space formulation studied here has received substantial attention: a sample of recent works, in addition to those cited below, include Wahba (1990), Johnstone and Silverman (1990), Koo (1993), Belitser and Levit (1995), Donoho (1995), Mair and Ruymgaart (1996), Golubev and Khas'minskiĭ (1999, 2001) and Cavalier, Golubev, Picard and Tsybakov (2002). However, much of this literature is concerned with eigenvalue sequences having (up to constants) monotonic behavior as $k$ increases. Papers that do specifically address the boxcar deconvolution problem include Hall, Ruymgaart, van Gaans and van Rooij (2001), Groeneboom and Jongbloed (2003) and O'Sullivan and Roy Choudhury (2001); see Section 5.1 for some further discussion.

1.2. *Effective degree of ill-posedness.* Problem (1) is an example of a linear statistical inverse problem in which one observes a noisy version of $Kf$ for some linear operator $K$, and wishes to reconstruct $f$. Such linear inverse problems are typically *ill-posed* in the sense of Hadamard: the inversion does not depend continuously on the observed data. One manifestation of this is that rates of convergence of estimators as $\epsilon \to 0$ are slower than in the direct case in which $f$ itself is observed with noise. We shall formulate some well-known existing results in terms of a notion of "degree of ill-posedness" (DIP) in order more easily to state the results of the present paper.

Under appropriate conditions, $K$ will have a singular value decomposition, and in terms of coefficients in the singular system expansions, the observations may be written in a sequence form

$$(9) \qquad y_k = r_k \theta_k + \epsilon z_k, \qquad k \in \mathbb{Z},$$

or, equivalently, after dividing through by $r_k$, as

$$(10) \qquad \bar{y}_k = \theta_k + \epsilon_k z_k,$$

where $\bar{y}_k = y_k/r_k$ and $\epsilon_k = \epsilon/r_k$. Let $\|\theta\|_2^2 = \sum_{k \in \mathbb{Z}} \theta_k^2$. Define the (nonlinear) minimax risk of estimation with respect to a parameter space $\Theta \subset \ell_2$ via

$$(11) \qquad R_N(\Theta, \epsilon) = \inf_{\hat{\theta}} \sup_{\theta \in \Theta} E \|\hat{\theta} - \theta\|_2^2,$$

where the infimum is taken over all (measurable) functions $\hat{\theta}$ of the data. We define the linear minimax risk by

$$R_L(\Theta, \epsilon) = \inf_{\hat{\theta}_L} \sup_{\theta \in \Theta} E \|\hat{\theta} - \theta\|_2^2,$$

where attention is restricted to the subclass of *linear* estimators $\hat{\theta}_L = (\hat{\theta}_k^L)$ with $\hat{\theta}_k^L = c_k y_k$, for some sequence $(c_k)$.

Parameter spaces of primary interest in this paper include, for $\sigma > 0, C > 0$, *hyperrectangles*

$$(12) \qquad H^\sigma(C) = \{\theta : |\theta_k| \leq C|k|^{-\sigma-1/2}, k \neq 0, \text{ and } \theta_0 \in \mathbb{R}\}$$

and *ellipsoids*

$$(13) \qquad \Theta_2^\sigma(C) = \left\{ \theta : \sum_k k^{2\sigma} \theta_k^2 \leq C^2 \right\}.$$

REMARK 1. Within these scales of spaces, the parameter $\sigma$ measures smoothness: larger $\sigma$ corresponds to faster decay of coefficients. When the $(\theta_k)$ are Fourier coefficients, the ellipsoids correspond exactly to mean-square smoothness of the $\sigma$ derivatives of $f = \sum \theta_k e_k$. [See, e.g., Kress (1999), Chapter 8.1.] There is no such simple characterization for hyperrectangles—the definition (12) is chosen to yield the same rates of convergence as (13) in the homogeneous cases described next. The parameter $C$ measures size: it corresponds to the radius of balls within these spaces.

REMARK 2. In (5) we used the complex exponentials $e^{\pi i k t}$. The model has the same form if instead one uses the real trigonometric basis $\bar{e}_k(t) = \cos \pi k t$ or $\sin \pi k t$ or $1/\sqrt{2}$ according as $k > 0, k < 0$ or $k = 0$. Model (9)–(10) applies to indices $k \in \mathbb{Z}$. For convenience in the rest of the paper, we restrict the index $k$ to $\mathbb{N}_+ = \{1, 2, \dots\}$. Indeed, since spaces such as (12) and (13) are symmetric with respect to $\pm k$, we have $R_N(\Theta, \epsilon; \mathbb{Z}) = 2R_N(\Theta, \epsilon; \mathbb{N}_+) + \epsilon^2$, with the analogous statement valid also for the linear minimax risks. Consequently, rates of convergence are certainly unaffected by working on $\mathbb{N}_+$.

REMARK 3. The notation $a(\epsilon) \asymp b(\epsilon)$ means that there exist constants such that for sufficiently small $\epsilon$, $c_1 b(\epsilon) \leq a(\epsilon) \leq c_2 b(\epsilon)$. The constants $c_1, c_2$ and other generic constants (denoted by $c$ and not necessarily the same at each appearance) may depend on parameters of the smoothness class $\Theta$ such as $\sigma$, but they do not depend on $\epsilon, \theta$ or the size parameter $C$. While the size constant $C$ clearly does not affect the rate of convergence as $\epsilon \to 0$, we consider it useful to show the order of dependence of minimax risks on $C$. The notation $a_k \sim b_k$ means that $\lim_{k \to +\infty}(a_k/b_k) = 1$. The notation $a_k \equiv c$ means that, for all $k$, $a_k = c$.

Suppose that the eigenvalues satisfy a homogeneous decay condition $r_k \sim |k|^{-\alpha}$ and that $\Theta = H^\sigma(C)$ or $\Theta = \Theta_2^\sigma(C)$. Then it is well known [e.g., Korostelev and Tsybakov (1993), Chapter 9] that

$$(14) \qquad R_N \asymp R_L \asymp C^{2(1-s)} \epsilon^{2s}, \qquad s = \frac{\sigma}{\sigma + 1/2 + \alpha}.$$

For direct data we have $r_k \equiv 1$ in (9) and it is known that

$$R_N \asymp R_L \asymp C^{2(1-s_D)} \epsilon^{2s_D}, \qquad s_D = \frac{\sigma}{\sigma + 1/2}.$$

This motivates the following definition of effective DIP:

$$(15) \qquad \alpha(K, \Theta) := \sigma\left(\frac{1}{s} - \frac{1}{s_D}\right).$$

For indirect problems $\alpha(K, \Theta)$ gives a measure of the effect (on the convergence rate) due to the inversion process. For example, if $K$ is an $\alpha$-fractional integration operator and $\Theta = \Theta_2^\sigma(C)$, then $r_k \sim |k|^{-\alpha}$ and so, in this case, $\alpha(K, \Theta) = \alpha$. As $\alpha$ gets larger it becomes more and more difficult to recover $f$.

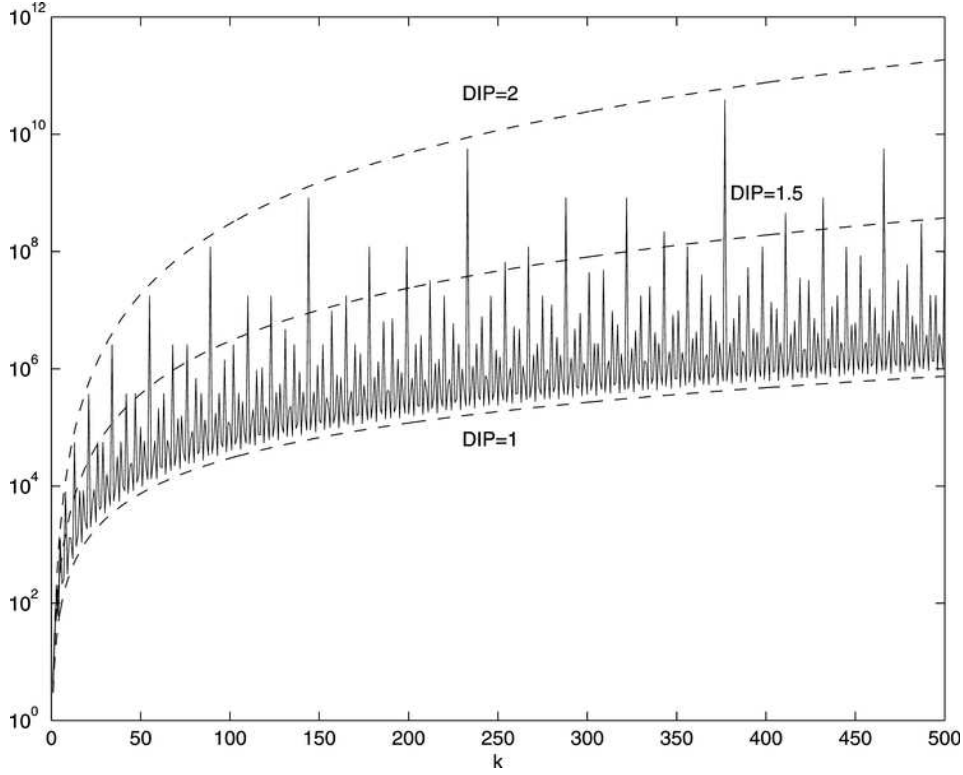

Fig. 1. *An illustration of the degree of the DIP for the boxcar deconvolution operator with $a = 2/(\sqrt{5}+1)$. Using a log-scale along the vertical axis, the function $k \to r_k^{-2}$ is depicted for $k = 0, 1, 2, \ldots, 500$ (oscillating solid line). For comparison purpose we also depict $k \to r_k^{-2}$ for a homogeneous operator with DIP= $1, 1.5, 2$ taking eigenvalues $r_k = ck^{-\alpha}$, where $\alpha = 1, 1.5, 2$ and $c = 0.58$ (smooth dashed curves).*

Returning to boxcar deconvolution, we note that $r_k \sim |k|^{-1}$ corresponds to an effective DIP of $\alpha = 1$. The question studied in this paper is whether the oscillations in $r_k$ of (3) increase the DIP. Compare Figure 1.

The answer turns out to depend on the function class. The main results, Theorems 1 and 2, can be expressed as saying, so long as logarithmic terms are ignored, that for ellipsoids and almost all irrational $a$,

$$\alpha(K_a, \Theta_2^\sigma) = \tfrac{3}{2} \qquad \text{for all } \sigma > 0,$$

while for hyperrectangles,

$$(16) \qquad \alpha(K_a, H^\sigma) = \begin{cases} 1, & \text{if } 0 < \sigma \leq \dfrac{3}{2}, \\ 1 + \dfrac{\sigma - 3/2}{2\sigma + 1}, & \text{if } \dfrac{3}{2} \leq \sigma. \end{cases}$$

Thus, the DIP of boxcar deconvolution lies between 1 and $\frac{3}{2}$, and is better (i.e., smaller) for more uniform smoothness (hyperrectangles) and for smaller $\sigma$.

REMARK 4. We caution that the literature contains other definitions of DIP of an inverse problem: for example, in Mathé and Pereverzev (2001), it refers to a numerical index of distance from invertibility. While these notions are certainly related, the definition used here is simply a convenience for interpreting results stated formally in Sections 3 and 4: it refers to the drop in rate of convergence due to presence of the decaying eigenvalues $r_k$.

REMARK 5. There is an elbow in rates at $\sigma = \frac{3}{2}$ for hyperrectangles but not ellipsoids. This contrasts with results obtained for homogeneous operators (14). Observe that the rates of convergence are worse for ellipsoids than for corresponding hyperrectangles: this occurs because the uniform hyperrectangle constraint (12) operates on *each* coordinate and so provides less scope for maximizing risk by concentrating signal energy in coordinates where $\|ka\|$ is small than does the ellipsoid case where only a total energy constraint (13) applies.

## 2. Preliminaries.

2.1. *Diophantine approximations.* We recall some pertinent parts of the classical theory, referring to Lang (1966) and Khinchin (1992) for further details. The study of approximations such as (8) is connected to the approximation of irrationals by rationals known as Diophantine approximations. For a given irrational number $a$, we distinguish the systematic approximations $\|ka\|$, $k = 1, 2, \ldots$ of (8) from the *best* rational approximations $p/q$: by *best*-approximation we mean that

$$(17) \qquad |qa - p| < \min_{1 \le k < q} \|ka\|.$$

Given the sequence of solutions $(p_n, q_n)$ to (17), the rate of approximation is defined in terms of the decay of

$$(18) \qquad D(a, q_n) = \left| a - \frac{p_n}{q_n} \right|.$$

Apart from the two basic groups of real numbers, rationals and irrationals, there exists a much finer division of irrational numbers based upon the degree to which they can be approximated by rational fractions. This may range from $O(1/q_n^2)$ to arbitrarily much faster, as explained below. These rates depend crucially on the best-possible rational approximation (17). The solution of (17) is given by the continued fractions algorithm which, unlike systematic fractions ($\|ka\|/k$, $k = 1, 2, \ldots$), captures the arithmetic properties of the number to be approximated.

2.2. *Continued fractions and convergents.* Any real number $a$ that is not an integer may be uniquely determined by its continued fraction expansion

$$(19) \qquad a = a_0 + \cfrac{1}{a_1 + \cfrac{1}{a_2 + \cfrac{1}{a_3 + \cdots}}} = [a_0; a_1, a_2, \ldots],$$

where $a_0$ is an integer and $a_1, a_2, \ldots$ is an infinite sequence of strictly positive integers. In the algorithm (19) the numbers $a_k$ are called the *elements* or *partial denominators*. To each infinite sequence $(a_k)$ corresponds a unique irrational number $a$ and vice versa. At stage $n$ the algorithm uses only the first $n$-*elements*: $[a_0; a_1, a_2, \ldots, a_n]$. For such a terminating continued fraction only a finite number of operations are involved and the result is clearly a rational number:

$$(20) \qquad a_0 + \cfrac{1}{a_1 + \cfrac{1}{\ddots + \cfrac{1}{a_n}}} = [a_0; a_1, a_2, \ldots, a_n] = \frac{p_n}{q_n}.$$

The rational numbers $(p_n/q_n)$, $n = 0, 1, \ldots$ are called the *convergents* of $a$. Returning to the problem of approximating an irrational number $a$ by rationals, we have that, for $n \geq 1$,

$$(21) \qquad \inf_{1 \leq k \leq q_n} \|ka\| = |q_n a - p_n| = \|q_n a\|.$$

In words, the *convergents* satisfy the best-approximation property (17). Indeed, any best-approximation is a convergent since, for $n \geq 1$, $q_n$ is the smallest integer $q > q_{n-1}$ such that $\|qa\| < \|q_{n-1}a\|$ [see, e.g., Lang (1966), page 9]. The quality of best-approximation is given by

$$(22) \qquad \frac{1}{2\,q_{n+1}} < \|q_n a\| < \frac{1}{q_{n+1}}$$

[Lang (1966), page 8]. While for systematic approximation, with $1 \leq k < q_n$, Lang [(1966), page 10] shows that

$$(23) \qquad \|ka\| > \frac{1}{2q_n}.$$

It is informative to note that, for $n \geq 2$, the algorithm (20) can be written as

$$(24) \qquad q_n = a_n q_{n-1} + q_{n-2}, \qquad p_n = a_n p_{n-1} + p_{n-2},$$

from which follow some basic properties of the convergents of all irrational numbers $a$:

(i) The denominators $q_n$ grow at least geometrically:

$$(25) \qquad q_{n+i} \geq 2^{(i-1)/2} q_i, \qquad i > 1.$$

(ii) For all $n \geq 0$,

$$a_n < \frac{q_n}{q_{n-1}} \leq a_n + 1.$$

The qualitative nature of rational approximations can, therefore, be measured by the size of the elements in the continued fraction algorithm, from (22),

(26) $$\frac{1}{2q_n^2(a_{n+1}+1)} < D(a, q_n) < \frac{1}{q_n^2 a_{n+1}}.$$

Faster approximation will occur for those irrationals with larger elements $a_n$ and vice versa. Families of irrational numbers can be defined according to the size of their elements.

DEFINITION 1. We say that an irrational number $a$ is badly approximable (BA) if

$$\sup_n a_n(a) < \infty.$$

From (26), we see that arbitrarily fast rates of approximation are possible.

A natural question arises—are there general laws which govern the approximations of classical irrational numbers?—Again, some answers follow from the continued fraction algorithm [Khinchin (1992), Chapter II]. One class of results concerns algebraic numbers—roots of polynomials with integer coefficients. For example, it can be shown that quadratic irrationals (such as $\sqrt{5}$) have periodic elements and so are BA. And cubic irrationals (e.g., $5^{1/3}$) cannot be approximated with a rate faster than $1/q^3$.

Another class of results constitutes the "measure theory" of continued fractions. For example, *almost all* numbers (i.e., except a set of Lebesgue measure zero) have unbounded $a_n$ [Khinchin (1992), Theorem 30]. On the other hand, for almost all numbers, it is also true that the rate of approximation can be no faster than $O(1/q_n^2 (\log q_n)^{1+\delta})$, $\delta > 0$. For us, an important consequence (see the Appendix) is the following. For each $\delta > 0$, there is a set $A_\delta$ of full measure such that

(27) $$q_{n+1} \geq q_n \log q_n \qquad \text{infinitely often,}$$

and yet

(28) $$q_{n+1} \leq q_n (\log q_n)^{1+\delta} \qquad \text{for all large } n > n(a).$$

Henceforth, the usage "almost all $a$" means "for all $a$ in $A_\delta$."

2.3. *Minimax risk.* We recall some basic results, established for the direct data setting $r_k \equiv 1$ (or $\epsilon_k \equiv \epsilon$) in Donoho, Liu and MacGibbon (1990), and easily extended to the indirect setting (10) (see the Appendix). If $\Theta$ is compact, orthosymmetric and quadratically convex, then

$$\text{(29)} \qquad R_N(\Theta, \epsilon) \leq R_L(\Theta, \epsilon) \leq \mu^\star R_N(\Theta, \epsilon),$$

where $\mu^\star \leq 1.25$ is the Ibragimov–Khasminskii constant; see Donoho, Liu and MacGibbon (1990). For such sets, we also have

$$\tfrac{1}{2} R_P(\Theta, \epsilon) \leq R_L(\Theta, \epsilon) \leq R_P(\Theta, \epsilon),$$

where we define

$$\text{(30)} \qquad R_P(\Theta, \epsilon) = \sup_{\theta \in \Theta} \sum_k \theta_k^2 \wedge \epsilon_k^2.$$

In the light of bounds (7) and Remark 2, our task is, then, to evaluate $R_P(\Theta, \epsilon)$ for selected $\Theta$, small $\epsilon$ and $k \in \mathbb{N}_+$, for the boxcar operator, which has

$$\text{(31)} \qquad \frac{\epsilon k}{\|ka\|} \leq \epsilon_k \leq \frac{\pi}{2} \frac{\epsilon k}{\|ka\|} \qquad \text{for all } k > 0.$$

2.4. *An equidistribution lemma.* While precise bounds (22) are available for best-possible rational approximations to an irrational number $a$, the quality of systematic rational approximations $\|ka\|$, $k = 1, 2, \ldots$, changes considerably as $k$ varies. As a result, $r_k$ and $r_k^{-2}$ oscillate widely as $k$ changes; see Figure 1. However, the *average* behavior is much less susceptible to fluctuations. Indeed, as $k$ runs over a block of length $q$, the values of $\|ka\|$ have a distribution that is in certain respects close to discrete uniform on $q^{-1}, 2q^{-1}, \ldots, 1$.

LEMMA 1. *Let $p/q$ and $p'/q'$ be successive principal convergents in the continued fraction expansion of a real number $a$. Let $N$ be a positive integer with $N + q < q'$. Let $h$ be a nonincreasing function. Then we have upper and lower bounds*

$$\text{(32)} \qquad \sum_{\mu=4}^{q} h(\mu/q) \leq \sum_{k=N+1}^{N+q} h(\|ka\|) \leq 2 \sum_{\mu=1}^{q-3} h(\mu/q) + 6h(1/(2q')).$$

PROOF. The argument is a modification of that used by [Lang (1966), page 37]. Since $p/q$ is a principal convergent, we may write $a$ in the form $a = p/q + \delta/q^2$ with $|\delta| < 1$. Writing $k = N + \nu$ with $\nu = 1, \ldots, q$, one gets

$$ka = Na + \nu p/q + \epsilon_\nu, \qquad |\epsilon_\nu| < 1/q.$$

Since $p$ and $q$ are relatively prime, the sets $\{\nu p/q, \nu = 1,\ldots,q\}$ and $\{\mu/q, \mu = 0,\ldots,q-1\}$ are equal modulo $\mathbb{Z}$. To each $k$ there is associated a unique $\nu$ and, hence, $\mu$, and setting $x_\mu = Na + \mu/q$, we have

$$ka = x_{\mu(k)} + \epsilon_{\mu(k)} \qquad (\mathrm{mod}\,\mathbb{Z}).$$

The points $\{x_\mu, \mu = 1,\ldots,q\}$ form an equispaced set with exactly one point in each interval $I_{\mu-1} = [(\mu-1)/q, \mu/q)$.

Let $R(\xi) = \xi - [\xi]$ denote the *remainder* of a real number $\xi$. Consider first the set $\mathcal{K}_1$ of indices $k$ for which the corresponding points $x_\mu$ lie in $I_0 \cup I_1 \cup I_{q-1}$: clearly, $|\mathcal{K}_1| = 3$. Since $k < q'$, we have from the remark following (22) that $R(ka) \geq \|ka\| \geq 1/(2q')$. Hence, the sum of $h(R(ka))$, for $k \in \mathcal{K}_1$, is bounded by $3h(1/(2q'))$.

Let $\mathcal{K}_2$ be the set of remaining indices $k$ in $\{N+1,\ldots,N+q\}$, so that the corresponding points $x_\mu$ lie in $I_2 \cup \cdots \cup I_{q-2}$. Since all $|\epsilon_\mu| < 1/q$, each of the left endpoints of $I_1,\ldots,I_{q-3}$ is a lower bound for exactly one $R(ka)$, $k \in \mathcal{K}_2$ and the right endoints of $I_3,\ldots,I_{q-1}$ each are upper bounds for exactly one $R(ka)$.

Combining this with the upper bound for $\mathcal{K}_1$, we obtain

$$(33) \qquad \sum_{\mu=4}^{q} h(\mu/q) \leq \sum_{k=N+1}^{N+q} h(R(ka)) \leq \sum_{\mu=1}^{q-3} h(\mu/q) + 3h(1/(2q')).$$

This inequality remains valid if we replace $h(R(ka))$ by $h(1 - R(ka))$—indeed, the proof is simply "reflected about $\frac{1}{2}$," and we note that for $k$ in the (reflected) $\mathcal{K}_1$, we have $1 - R(ka) \geq \|ka\| > 1/(2q')$. Since $\|x\| = \min\{R(x), 1 - R(x)\}$, we have

$$h(\|x\|) = \max\{h(R(x)), h(1 - R(x))\},$$

and using $(a+b)/2 \leq \max\{a,b\} \leq a+b$, the lemma follows from (33) applied to $R(ka)$ and $1 - R(ka)$. $\square$

REMARK 6. The proof shows that the upper bound continues to hold if the middle sum is taken over $N + 1 \leq k \leq N + k_0$, where $k_0 \leq q$ and we assume only $N + k_0 < q'$.

REMARK 7. The bounds provided by this lemma are often sharp up to constants. For example, if $a$ is BA and $h(x) = 1/x$,

$$\sum_{k=N+1}^{N+q} \|ka\|^{-1} \asymp q \log q.$$

**3. Hyperrectangles.**

3.1. *Statement and outline.* To state the main results, introduce two rate constants

$$r = (\sigma + \tfrac{1}{2})/(\sigma + \tfrac{5}{2}), \qquad \bar{r} = \sigma/(\sigma + \tfrac{3}{2}),$$

and note that $r < \bar{r}$ if and only if $\sigma > 3/2$. More precise results are possible in the BA case, while for generic irrationals, the consequences (27) and (28) of Khinchin's theorem lead to only slightly weaker statements.

THEOREM 1. *For BA a we have*

(34) $$R_N(H^\sigma(C), \epsilon) \asymp \begin{cases} C^{2(1-r)} \epsilon^{2r}, & \text{if } \sigma > \tfrac{3}{2}, \\ C\epsilon \log(C/\epsilon), & \text{if } \sigma = \tfrac{3}{2}, \\ C^{2(1-\bar{r})} \epsilon^{2\bar{r}}, & \text{if } 0 < \sigma < \tfrac{3}{2}. \end{cases}$$

*For almost all a, the previous bounds remain valid for $0 < \sigma < \tfrac{3}{2}$, while for $\sigma \geq \tfrac{3}{2}$, for each $\delta > 0$,*

(35) $R_N(H^\sigma(C), \epsilon) \begin{cases} \leq c_2 (\log C/\epsilon)^{5+\delta} C^{2(1-r)} \epsilon^{2r} & \text{for all small } \epsilon, \\ \geq c_1 (\log C/\epsilon)^{2r} C^{2(1-r)} \epsilon^{2r} & \text{for infinitely many } \epsilon. \end{cases}$

There is thus an "elbow" in the rates of convergence at $\sigma = \tfrac{3}{2}$. Comparison with (14) shows that for $\sigma < \tfrac{3}{2}$, the DIP is $\alpha = 1$ (as if the sinusoidal term were not present in $r_k$). However, for $\sigma > \tfrac{3}{2}$, the DIP given by (16) increases gradually from 1 to a limiting value of $\tfrac{3}{2}$ for large $\sigma$.

This result does not cover irrationals with fast rates of approximation (e.g., $1/q^3$ or higher, as discussed in Section 2.2), but, of course, such numbers form a set of Lebesgue measure zero.

We outline the main steps of the proof, with details to follow in Section 3.3. First, as notational convention, we introduce a parameter $\tau = \sigma + \tfrac{1}{2}$, so that $\Theta = H^{\tau - 1/2}(C) = \{\theta : |\theta_k| \leq Ck^{-\tau}\}$. With these conventions, (30) becomes

(36) $$R_P(\Theta, \epsilon) = \sum_{k>0} C^2 k^{-2\tau} \wedge \epsilon_k^2 := \sum_{k>0} m_k(\epsilon).$$

First, we use the continued fraction approximation to $a$: $p_n/q_n, n = 0, 1, 2, \ldots$, and for frequencies near $q_n$, split the sum into blocks of length $q_n$. Thus,

(37) $$\sum_{k>0} m_k(\epsilon) = \sum_{blocks} \sum_{k \in block} m_k(\epsilon),$$

where $\sum_{blocks}$ is the sum over all blocks as $n$ varies, the blocks being of length $q_n$ between $q_n$ and $q_{n+1}$. We then apply the equidistribution lemma to the sum within blocks. The block sums are then collected into one of three zones:

(38) $$R_P(\Theta, \epsilon) = \sum_k m_k(\epsilon) = V(\epsilon) + M(\epsilon) + B(\epsilon).$$

These zones (variance, mixed and bias) are illustrated in Figure 2, and defined formally at (45).

3.2. *Frequency partitions determined by an irrational.* Any irrational number $a$ defines a unique sequence of convergents: $p_n/q_n; 1 = q_0 < q_1 < \cdots < q_n < q_{n+1} < \cdots$. Define $l_n \geq 1$ as the largest integer strictly less than $q_{n+1}/q_n$, thus,

$$l_n q_n < q_{n+1} \leq (l_n + 1) q_n.$$

Consider a nonuniform grid

$$\ldots, q_n, 2q_n, \ldots, l_n q_n, \quad q_{n+1}, 2q_{n+1}, \ldots, l_{n+1} q_{n+1}, \quad q_{n+2}, \ldots.$$

Introduce indices $\nu = (n, l), l = 1, \ldots, l_n; n = 1, 2, \ldots$. The bivariate indices $\nu = (n, l)$ are totally ordered by lexicographic ordering and we refer to their components by the functions $n(\nu)$, $l(\nu)$. Furthermore, each index $\nu$ has an immediate successor, which in slight abuse of notation we denote by $\nu + 1$. So our grid is

(39) $$N_\nu = l(\nu) q_{n(\nu)};$$

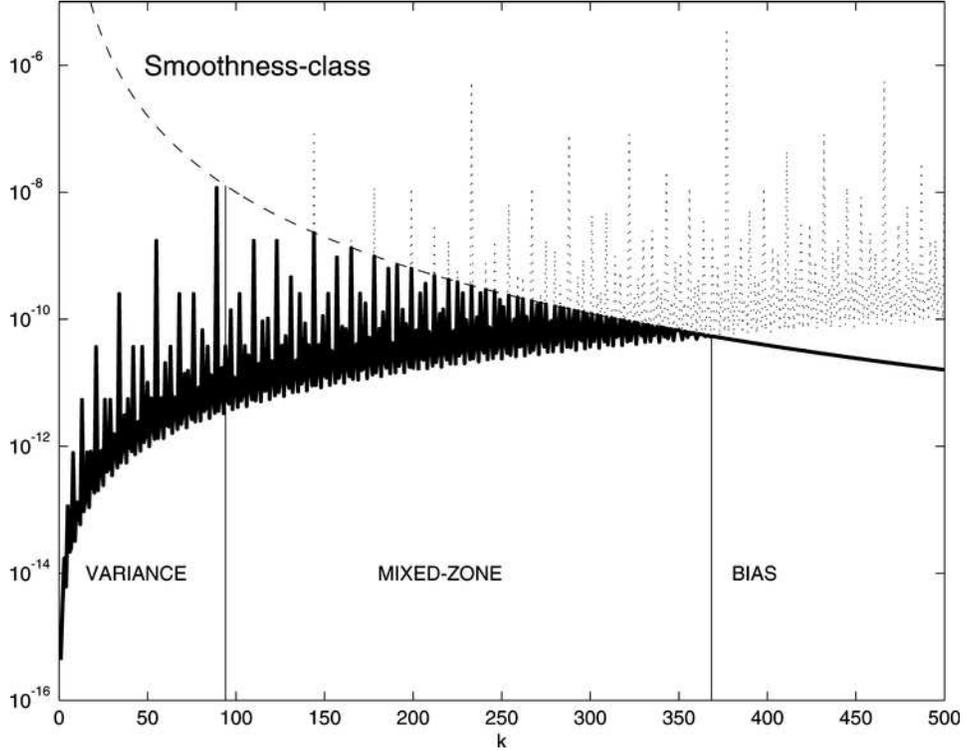

FIG. 2. *An illustration of the variance-mixed-bias zones. Using a log-scale along the vertical axis, the plot shows both functions $k \to \epsilon_k^2$ (oscillating dotted curve) and $k \to C^2 k^{-2\tau}$ (smooth dashed curve), with $a = 2/(\sqrt{5}+1)$, $\epsilon = 10^{-8}$, $C = 1$ and $\tau = 2$, which corresponds to $\sigma = 3/2$. Solid vertical lines indicate the borders of the key zones. The thick solid line plots $k \to m_k(\epsilon) = C^2 k^{-2\tau} \wedge \epsilon_k^2$.*

this grid defines a partition of $\mathbb{N}_+$ by blocks which between $q_n$ and $q_{n+1}$ have length $\leq q_n$:

$$(40) \qquad \mathbb{N}_+ = \bigcup_\nu B_\nu, \qquad B_\nu = [N_\nu, N_{\nu+1}).$$

Clearly,

$$|B_\nu| = N_{\nu+1} - N_\nu = \begin{cases} q_{n(\nu)}, & \text{unless } l(\nu) = l_{n(\nu)}, \\ \in [1, q_{n(\nu)}), & \text{if } l(\nu) = l_{n(\nu)}. \end{cases}$$

To simplify certain calculations we use blocks of length $q_{n(\nu)}$ only, introducing

$$(41) \qquad C_\nu = [N_\nu, N_\nu + q_{n(\nu)}] \supset B_\nu.$$

By construction, for a given integer $k$, there are at most two $C_\nu$ such that $k \in C_\nu$. Hence, summing over all $C_\nu$ in place of $B_\nu$ will only affect the rate by a multiplicative constant of at most 2.

### 3.3. Proof of Theorem 1.

3.3.1. *Key zones and bounds.* First, recall that $m_k(\epsilon)$ is defined at (36) and use bounds (31); by construction $q_{n(\nu)} \leq N_\nu$ so that for $k$ in a block $[N_\nu, N_\nu + q_{n(\nu)}]$, $N_\nu \leq k \leq 2N_\nu$, hence,

$$(42) \quad m_k(\epsilon) \asymp C^2 k^{-2\tau} \wedge \epsilon^2 \frac{k^2}{\|ka\|^2} \asymp C^2 N^{-2\tau} \wedge \epsilon^2 \frac{N^2}{\|ka\|^2} := h_N(\|ka\|).$$

We suppress the index $\nu$ when not necessary. From the equidistribution lemma,

$$(43) \qquad \sum_{\mu=4}^{q} h_N\left(\frac{\mu}{q}\right) \leq \sum_{k \in C_\nu} h_N(\|ka\|) \leq c \sum_{\mu=1}^{q} h_N\left(\frac{\mu}{q}\right) + c h_N\left(\frac{1}{2q'}\right).$$

To estimate these sums, we use an easily verified bound.

LEMMA 2. *If $q > 2r$ and $\kappa > 0$, then*

$$\sum_{\mu=r}^{q} 1 \wedge \left(\frac{\kappa}{\mu}\right)^2 \asymp \min\{\kappa^2, \kappa, q\},$$

*where the constants needed for $\asymp$ depend only on $r$.*

Now apply this to $h_N(x) = C^2 N^{-2\tau} \wedge \epsilon^2 N^2 x^{-2}$. Writing also $\varepsilon = \epsilon/C$, we obtain

$$(44) \qquad \sum_{\mu=r}^{q} h_N(\mu/q) \asymp C^2 N^{-2\tau} \min\{\varepsilon^2 N^{2(1+\tau)} q^2, \varepsilon N^{1+\tau} q, q\}.$$

We can now formally define the zone to which a block $B_\nu$ (or $C_\nu$) belongs in terms of the value of $\varepsilon N_\nu^{1+\tau} q_{n(\nu)}$. Again suppressing the subscript $\nu$, we say

$$B_\nu \in \begin{cases} \text{Variance zone} & \Leftrightarrow \quad \varepsilon N^{1+\tau} q \le 1, \\ \text{Mixed zone} & \Leftrightarrow \quad 1 < \varepsilon N^{1+\tau} q \le q, \\ \text{Bias zone} & \Leftrightarrow \quad \varepsilon N^{1+\tau} q > q. \end{cases} \quad (45)$$

Thus, the zone describes which term appears in the minimizer in (44). Let $\nu_0 < \nu_1$ be the last indices for which $\varepsilon N_\nu^{1+\tau} q_{n(\nu)} \le 1$ and $\varepsilon N_\nu^{1+\tau} \le 1$, respectively, and set

$$k_0(\varepsilon) = N_{\nu_0 + 1} \quad \text{and} \quad k_1(\varepsilon) = N_{\nu_1 + 1}. \quad (46)$$

Frequencies $k < k_0$ lie in the variance zone, those with $k_0 \le k < k_1$ in the mixed zone, and those with $k \ge k_1$ in the bias zone.

Consider now the second term in the upper bound of (43):

$$h_N(1/(2q')) = C^2 N^{-2\tau}(1 \wedge (2\varepsilon N^{1+\tau} q')^2).$$

If $\varepsilon N^{1+\tau} q > 1$, then, of course, so is $\varepsilon N^{1+\tau} q'$ and so $h_N(1/(2q')) = C^2 N^{-2\tau} < C^2 \varepsilon N^{1-\tau} q$ can be ignored in comparison with (44). On the other hand, if $\varepsilon N^{1+\tau} q \le 1$, then $h_N(1/(2q')) \le 4\epsilon^2 N^2 (q')^2$ and this bound dominates $\epsilon^2 N^2 q^2$. In summary, we have derived the following key bounds:

$$\sum_{k \in C_\nu} m_k(\epsilon) \begin{cases} \le c\epsilon^2 N^2 (q')^2, & \nu \in \text{(variance zone)}, \\ \asymp C\epsilon N^{1-\tau} q, & \nu \in \text{(mixed zone)}, \\ \asymp C^2 N^{-2\tau} q, & \nu \in \text{(bias zone)}. \end{cases} \quad (47)$$

*The variance zone.* Consider first values $k < k_0(\epsilon)$ such that the contribution to the minimax risk is due to oscillations occasioned by Diophantine approximation only. Here the first bound of (47) applies and the hyperrectangle constraint $k \to C^2 k^{-2\tau}$ has not yet any smoothing effect.

We first derive an expression for $k_0$ in terms of $\epsilon$. If $\nu = \nu_0 + 1$, we have by definition,

$$\varepsilon^{-1} < N_\nu^{1+\tau} q_{n(\nu)} \le N_\nu^{2+\tau} = k_0^{2+\tau} \quad \text{and so } k_0 \ge \varepsilon^{-1/(2+\tau)}.$$

On the other hand, again by definition, $\varepsilon^{-1} > N_{\nu_0}^{1+\tau} q_{n(\nu_0)} \ge q_{n(\nu_0)}^{2+\tau}$. Writing $L_n$ for $q_{n+1}/q_n$, we obtain

$$k_0 = N_{\nu_0 + 1} \le q_{n(\nu_0)+1} \le L_{n(\nu_0)} q_{n(\nu_0)} \le L_{n(\nu_0)} \varepsilon^{-1/(2+\tau)}.$$

For BA $a$, $L_n \le c$, while for almost all $a$ and all large $n$, (28) shows that $L_n \le (\log q_n)^{1+\delta}$. To summarize,

$$k_0 \asymp (C/\epsilon)^{1/(2+\tau)} \quad \text{for BA } a,$$
$$k_0 \le c(C/\epsilon)^{1/(2+\tau)}((\log C/\epsilon))^{1+\delta} \quad \text{for almost all } a.$$

First, sum over blocks using partition (40) and apply bound (47) in the variance zone:

$$(48) \quad V(\epsilon) = \sum_{k=1}^{k_0-1} m_k(\epsilon) = \sum_{\nu \leq \nu_0} \sum_{k \in B_\nu} m_k(\epsilon) \leq c\epsilon^2 \sum_{\nu \leq \nu_0} N_\nu^2 q_{n(\nu)+1}^2.$$

Using grid (39), and setting $\nu_0 = (n_0, l), l \leq l_{n_0}$, we obtain

$$(49) \quad V(\epsilon) \leq c\epsilon^2 \sum_{n=1}^{n_0} \sum_{l=1}^{l_n} l^2 q_n^2 q_{n+1}^2 \leq c\epsilon^2 \sum_{n=1}^{n_0} l_n^3 q_n^2 q_{n+1}^2 \leq c\epsilon^2 \bar{L}_{n_0}^5 \sum_{n=1}^{n_0} q_n^4,$$

where we have set $\bar{L}_{n_0} = \max\{L_n, n \leq n_0\}$.

The denominators $q_n$ grow at least exponentially [cf. (25)] and so using $q_{n_0} \leq \varepsilon^{-1/(2+\tau)}$, we find

$$\epsilon^2 \sum_{n=1}^{n_0} q_n^4 \leq c\epsilon^2 q_{n_0}^4 \leq c\epsilon^2 (C/\epsilon)^{4/(2+\tau)} = cC^{2(1-r)}\epsilon^{2r}.$$

In the BA case, $\bar{L}_{n_0} \leq c$, while for almost all $a$ we have $\bar{L}_{n_0} \leq (\log q_{n_0})^{1+\delta/5} \leq c(\log \varepsilon)^{1+\delta/5}$. In summary,

$$(50) \quad V(\epsilon) \leq \begin{cases} cC^{2(1-r)}\epsilon^{2r}, & \text{for BA } a, \\ c(\log(C/\epsilon))^{5+\delta} C^{2(1-r)}\epsilon^{2r}, & \text{for almost all } a. \end{cases}$$

*The mixed zone.* We are now interested in indices $k \in [k_0, k_1)$ where both oscillations and the hyperrectangle constraint $C^2 k^{-2\tau}$ contribute to the minimax risk; it ends where the oscillations stop. By definition, $k_1 = N_{\nu_1+1}$ satisfies $N_{\nu_1} \leq \varepsilon^{-1/(1+\tau)} < N_{\nu_1+1}$. Since always $N_{\nu+1} \leq 2N_\nu$, it follows that

$$k_1 \asymp \varepsilon^{-1/(1+\tau)} = (C/\epsilon)^{1/(1+\tau)}.$$

Using bound (47) in the mixed zone, together with $|C_\nu| = q_{n(\nu)}$, and $N \leq k \leq 2N$ yields

$$\sum_{k \in C_\nu} m_k(\epsilon) \asymp C\epsilon N_\nu^{1-\tau} q_{n(\nu)} \asymp C\epsilon \sum_{k \in C_\nu} N_\nu^{1-\tau} \asymp C\epsilon \sum_{k \in C_\nu} k^{1-\tau},$$

which shows that for sums over blocks of length $q_{n(\nu)}$ in the mixed zone, we may replace $m_k(\epsilon)$ by $\epsilon k^{1-\tau}$. Since the blocks $C_\nu$ form a cover of the integers $k_0, \ldots, k_1 - 1$ of redundancy at most two,

$$M(\epsilon) = \sum_{k=k_0}^{k_1-1} m_k(\epsilon) \asymp C\epsilon \sum_{k=k_0}^{k_1-1} k^{1-\tau}.$$

Thus, in the mixed zone,

$$(51) \quad M(\epsilon) \asymp \begin{cases} C\epsilon k_0^{2-\tau} \asymp C^{2(1-r)}\epsilon^{2r}, & \text{if } \tau > 2, \\ C\epsilon \log(k_1/k_0) \asymp C\epsilon \log(C/\epsilon), & \text{if } \tau = 2, \\ C\epsilon k_1^{2-\tau} \asymp C^{2(1-\bar{r})}\epsilon^{2\bar{r}}, & \text{if } \frac{1}{2} < \tau < 2. \end{cases}$$

*The bias zone.* Note that for $k \geq k_1$, since always $\|ka\| \leq 1$, we have $\epsilon^2 k^2 / \|ka\|^2 \geq \epsilon^2 k^2 \geq C^2 k^{-2\tau}$ and so there is no longer any effect of oscillation, and $m_k(\epsilon) = C^2 k^{-2\tau}$ in (36). Hence,

$$(52) \qquad B(\epsilon) = \sum_{k \geq k_1} m_k(\epsilon) = C^2 \sum_{k \geq k_1} k^{-2\tau} \asymp C^2 k_1^{-2\tau+1} \asymp C^{2(1-\bar{r})} \epsilon^{2\bar{r}}.$$

We emphasize that bounds (51) and (52) apply to all irrationals $a$.

3.3.2. *Summary.* We return to (38). In the BA case (and also the a.a. case when $\frac{1}{2} < \tau < 2$), it is apparent from (50), (51) and (52) that $V + B + M \asymp M$, which establishes (34).

It remains to consider the a.a. case with $\tau \geq 2$. The upper bound in (35) is apparent from (50). For the lower bound, let $a$ be an arbitrary irrational with convergents $p_k/q_k$, $k = 0, 1, 2, \ldots$. Simply by choosing $\theta$ to be zero except in the $k$th coordinate—in which $\theta_k = Ck^{-\tau}$—we obtain the elementary lower bound

$$(53) \qquad R_P(\Theta, \epsilon) \geq \sup_k C^2 k^{-2\tau} \wedge \epsilon_k^2.$$

Since $\epsilon_k \geq \epsilon k/\|ka\|$, we find using (22) that for $k = q_n$,

$$C^2 k^{-2\tau} \wedge \epsilon_k^2 \geq C^2 q_n^{-2\tau} \wedge \epsilon^2 q_n^2 q_{n+1}^2.$$

Using (27) in (53), we deduce that for almost all $a$ there exists a sequence $n_l$ such that

$$(54) \qquad R_P(\Theta, \epsilon) \geq \sup_l C^2 q_{n_l}^{-2\tau} \wedge \epsilon^2 q_{n_l}^4 (\log q_{n_l})^2.$$

Construct a sequence $(\epsilon[l]), l = 1, 2, \ldots$, with

$$(55)\ C^2 q_{n_l}^{-2\tau} = \epsilon[l]^2 q_{n_l}^4 (\log q_{n_l})^2, \qquad \text{which gives } q_{n_l} \asymp (\varepsilon[l] \log \varepsilon[l]^{-1})^{-1/(2+\tau)},$$

and using such an $\epsilon[l]$-sequence in (54), together with (55), yields the required bound

$$R_P(\Theta, \epsilon[l]) \geq C^2 q_{n_l}^{-2\tau} \asymp (\log(C/\epsilon[l]))^{2r} C^{2(1-r)} \epsilon[l]^{2r}.$$

**4. Ellipsoids.** For an ellipsoid $\Theta = \Theta^\sigma(C)$ defined as in (13), let $\tilde{r} = \sigma/(\sigma + 2)$. The goal of this section is to establish the following:

THEOREM 2. *For $\sigma > 0$ and BA $a$, we have*

$$R_N(\Theta^\sigma(C), \epsilon) \asymp C^{2(1-\tilde{r})} \epsilon^{2\tilde{r}}.$$

*For almost all $a$, bounds (35) hold for $R_N(\Theta^\sigma(C), \epsilon)$ with $r$ replaced by $\tilde{r}$, for all $\sigma > 0$.*

Since $\tilde{r} = \sigma/(\sigma+2)$, the DIP $\alpha(K_a, \Theta^\sigma(C)) = \frac{3}{2}$ for all ellipsoids, regardless of the value of the smoothness index $\sigma$.

*Upper bound.* As with hyperrectangles, the aim is to use sums over blocks of length $\approx q$. To do so, we define slightly larger ellipsoids based on the partition $\{B_\nu\}$ of (40):

$$\Theta_a = \Theta_a^\sigma(C) = \left\{\theta : \sum_\nu N_\nu^{2\sigma} \sum_{k \in B_\nu} \theta_k^2 \leq C^2\right\}, \tag{56}$$

where the index $a$ indicates that the grid depends on number theoretical properties of $a$. By definition (40) of the partition, $k \in B_\nu$ implies that $k \geq N_\nu$ so that $\Theta \subset \Theta_a$ and, hence, $R(\Theta, \epsilon) \leq R(\Theta_a, \epsilon)$.

We may now split the optimization across and within blocks:

$$\begin{aligned}
R_P(\Theta_a, \epsilon) &= \sup\left\{\sum_\nu \sum_{k \in B_\nu} \theta_k^2 \wedge \epsilon_k^2 : \theta \in \Theta_a\right\} \\
&= \sup\left\{\sum_\nu b_\nu(t_\nu, \epsilon) : \sum_\nu N_\nu^{2\sigma} t_\nu^2 \leq C^2\right\},
\end{aligned} \tag{57}$$

where the optimization within block $B_\nu$ is subject to the quota $t_\nu^2$:

$$b_\nu(t_\nu, \epsilon) = \sup\left\{\sum_{k \in B_\nu} \theta_k^2 \wedge \epsilon_k^2 : \sum_{k \in B_\nu} \theta_k^2 \leq t_\nu^2\right\} = \min\left\{t_\nu^2, \sum_{k \in B_\nu} \epsilon_k^2\right\}. \tag{58}$$

The equidistribution lemma can be applied to this last sum: $\sum \epsilon_k^2 \asymp \epsilon^2 \sum k^2/\|ka\|^2$. On dropping the subscript $\nu$, we obtain

$$\sum_{k \in B_\nu} \frac{k^2}{\|ka\|^2} \leq 4N^2 \sum_{k \in B_\nu} \frac{1}{\|ka\|^2} \leq 8N^2 \left\{\sum_{\mu=1}^q \frac{q^2}{\mu^2} + 3(2q')^2\right\} \leq cN^2(q')^2.$$

Hence, from (57) and (58),

$$R_P(\Theta_a, \epsilon) \leq c \sup\left\{\sum_\nu \min\{t_\nu^2, \epsilon^2 N_\nu^2 q_{n(\nu)+1}^2\} : \sum_\nu N_\nu^{2\sigma} t_\nu^2 \leq C^2\right\}. \tag{59}$$

Observe that for any positive sequences $(u_\nu), (c_\nu)$ and $(d_\nu)$, with $d_\nu$ nondecreasing,

$$\sup_u \left(\sum_\nu \min(u_\nu, c_\nu) : \sum_\nu d_\nu u_\nu \leq 1\right) \leq \sum_{\nu \leq \nu_0} c_\nu \tag{60}$$

for any value $\nu_0$ for which

$$\sum_{\nu \leq \nu_0} c_\nu d_\nu \geq 1. \tag{61}$$

Applying this to (59) with $u_\nu = t_\nu^2, c_\nu = \epsilon^2 N_\nu^2 q_{n(\nu)+1}^2$ and $d_\nu = N_\nu^{2\sigma}/C^2$, we obtain

(62) $$R_P(\Theta_a, \epsilon) \leq c\epsilon^2 \sum_{\nu \leq \nu_0} N_\nu^2 q_{n(\nu)+1}^2.$$

Here $c_\nu d_\nu = \varepsilon^2 N_\nu^{2\sigma+2}(q')^2 = \varepsilon^2 (lq_n)^{2\sigma+2} q_{n+1}^2$ if $\nu = (n, l)$. Let $\mathcal{N}_n = \{\nu : q_n \leq N_\nu < q_{n+1}\}$ and note, since $(l_n + 1)q_n \geq q_{n+1}$, that

$$CD_n := \sum_{\nu \in \mathcal{N}_n} c_\nu d_\nu = \varepsilon^2 q_n^{2\sigma+2} q_{n+1}^2 \sum_1^{l_n} l^{2\sigma+2}$$

$$\geq c\varepsilon^2 (l_n + 1)^{2\sigma+3} q_n^{2\sigma+2} q_{n+1}^2 \geq c\varepsilon^2 l_n q_{n+1}^{2\sigma+4}.$$

Let $n_0$ be the first index $n$ for which $CD_n \geq 1$: since $CD_{n_0-1} \leq 1$, we have

(63) $$\varepsilon^2 q_{n_0}^{2\sigma+4} < 1/(cl_{n_0-1}) \quad \text{and so } q_{n_0} \leq c\varepsilon^{-1/(\sigma+2)}.$$

Since (62), together with (63), is exactly the situation reached at (48) in the hyperrrectangle case (with $\tau$ replaced by $\sigma$) we conclude that the bounds (50) apply (with $r$ replaced by $\tilde{r}$).

*Lower bound.* Arguing exactly as at (53), but with $\tau$ replaced by $\sigma$,

(64) $$R_P(\Theta, \epsilon) \geq \sup_n \ C^2 q_n^{-2\sigma} \wedge \epsilon^2 q_n^2 q_{n+1}^2.$$

In the BA case, let $n_0$ be the last index $n$ for which $\epsilon^2 q_n^4 < C^2 q_n^{-2\sigma}$, so that $q_{n_0}^{2\sigma+4} < \varepsilon^{-2}$ and $q_{n_0}^{-2\sigma} > \varepsilon^{2\sigma/(\sigma+2)}$. From (64) at $n = n_0 + 1$, we find

$$R_P(\Theta, \epsilon) \geq C^2 q_{n_0+1}^{-2\sigma} \geq cC^2 q_{n_0}^{-2\sigma} \geq cC^{2(1-\tilde{r})} \epsilon^{2\tilde{r}}.$$

For the almost all case, the argument is the same as before at (54) and below.

## 5. Discussion.

5.1. *Periodic vs. nonperiodic.* Recent papers by Hall, Ruymgaart, van Gaans and van Rooij (2001) and Groeneboom and Jongbloed (2003) consider in part a density estimation version of the deconvolution problem in which the data consist of an i.i.d. sample $Y_i = X_i + z_i$ in which $X_i$ are i.i.d. with unknown density $f$ and $z_i$ are i.i.d. uniform on $[-a, a]$ and independent of the $X_i$. Groeneboom and Jongbloed (2003) derive pointwise limiting distributions of estimators of $f$ based on kernel smooths of nonparametric MLEs of the distribution function of $f$. The work of Hall, Ruymgaart, van Gaans and van Rooij (2001) looks at maximum global estimation errors, and so is perhaps closer in spirit to the present investigation. Instead of any periodicity assumptions, it is assumed there that the density $f$ has compact support on $\mathbb{R}$. The compact support permits an explicit inversion formula: if $g = K_a f$ and $I$ is chosen large enough that $x - Ia < \inf \operatorname{supp} f$, then

$$f(x) = 2a \sum_{i=1}^{I} g'(x - ia).$$

In this case Hall, Ruymgaart, van Gaans and van Rooij (2001) show that the DIP $\alpha(K_a, \mathcal{F}^\sigma) = 1$ for $\mathcal{F}^\sigma$ of both hyperrectangle and ellipsoid type, in contrast to the results found for the periodic model considered here. The difference in results may perhaps be understood by observing that sinusoids, which are basic to the periodic model, do not have compact support. Thus, the models capture genuinely different phenomena.

5.2. *Effect of rational approximations to $a$.* In practice, computer code works with rational numbers—what effect will this have on our conclusions? A few remarks can be made even without getting into specifics of particular models of computation or attempting a full analysis.

A basic issue is whether the boxcar width $a$ is under the investigator's control. If it is—our first scenario—then we might imagine replacing $a$ by $\alpha_m = p_m/q_m$, say, so that model (4) becomes

$$(65) \qquad y_k = r_k(\alpha_m)\theta_k + \epsilon z_k, \qquad r_k(\alpha_m) = \frac{\sin \pi k \alpha_m}{\pi k \alpha_m}.$$

Here $p_m/q_m$ might be one of the sequence of best rational approximations to $a$. The approximation results of Section 2.2 show that our analysis of estimation in model (65) is unchanged from that of irrational $a$, at least for frequencies $k \leq q_m$, since $a$ and $\alpha_m$ will have the same convergents $p_r/q_r$ for $r \leq m$. Thus, one could simply choose $q_m$ large enough that the tail bias accruing to frequencies above $q_m$ is negligible. To be more specific, assume that $\Theta$ is a hyperrectangle $H^\sigma(C)$, and that $\epsilon$ is known. Let $\eta > 0$ be small [we could let $\eta(\epsilon) \to 0$ with $\epsilon$ to preserve rates of convergence]. We can choose $k_2 > k_1(\epsilon)$ [defined at (46)] so that the tail bias

$$C^2 \sum_{k > k_2} k^{-2\sigma+1} \leq \eta R(H^\sigma(C), \epsilon),$$

and then choose $m$ large enough that $q_m \geq k_2$. A minimax estimator for $H^\sigma(C)$ under model (65) will be essentially identical in structure with one for the original irrational $a$, since in either case, the zero estimator is used at all frequencies $k > q_m$.

In the second scenario, the boxcar width $a$ is determined by nature and the investigator must work with the data $y$ from model (4). We still assume that the value of $a$ is known, but must use rational approximations to $a$ in our estimators based on $y$. For definiteness, consider again the case $\Theta = H^\sigma(C)$ and set $\tau_k = Ck^{-\tau}$. Consider the risk of linear rules $\hat{\theta}_k(y) = c_k y_k$ if $\epsilon_k \leq \tau_k$ and $\hat{\theta}_k(y) = 0$ otherwise. If $\mathcal{S} = \{k : \epsilon_k \leq \tau_k\}$, then the risk of such a rule is

$$r(c, \theta) = \sum_{k \in \mathcal{S}} [c_k^2 \epsilon^2 + (1 - c_k r_k)^2 \theta_k^2] + \sum_{k \notin \mathcal{S}} \theta_k^2.$$

Suppose that $a$ is irrational: with infinite precision, we could use an estimator $c_k = 1/r_k$ that makes $r(c, \theta) = \sum \theta_k^2 \wedge \epsilon_k^2$. Now consider the difference in risk that results from an approximation $\hat{c}_k = 1/\hat{r}_k$, where $\hat{r}_k = (\sin \pi k \hat{a})/(\pi k \hat{a})$ for some rational approximation $\hat{a} = p_m/q_m$ to $a$,

$$r(\hat{c}, \theta) - r(c, \theta) = \sum_{\mathcal{S}} \left\{ \left[ \left( \frac{r_k}{\hat{r}_k} \right)^2 - 1 \right] \epsilon_k^2 + \left( 1 - \frac{r_k}{\hat{r}_k} \right)^2 \theta_k^2 \right\};$$

if we write $r_k/\hat{r}_k = 1 + \delta_k$, and assume that $\bar{\delta} = \sup_{k \in \mathcal{S}} |\delta_k| \leq 1$,

(66) $$\sup_\Theta |r(\hat{c}, \theta) - r(c, \theta)| \leq 3\bar{\delta} R_P(\Theta, \epsilon) + \bar{\delta}^2 \sum \tau_k^2.$$

Using a derivative bound on $a \to \sin \pi k a$ and then (7),

$$|\delta_k| \leq \frac{\hat{a}}{a} \left| \frac{\sin \pi k a}{\sin \pi k \hat{a}} - 1 \right| + \left| \frac{\hat{a}}{a} - 1 \right| \leq \frac{|\hat{a} - a|}{a} \left\{ \frac{\pi k}{\sin \pi k \hat{a}} + 1 \right\} \leq \frac{2|\hat{a} - a|}{a} \frac{k}{\|k\hat{a}\|}.$$

If $\hat{a} = p_m/q_m$ and $k < q_r$, then from (26), (23) and (25),

$$|\delta_k| \leq \frac{4}{a} \left( \frac{q_r}{q_m} \right)^2 \leq \frac{8}{a} 2^{-(m-r)}.$$

Consequently, the risk difference due to using a rational approximation $\hat{a}$ can be made as small as desired by first selecting $r$ so that $\sup\{k : k \in \mathcal{S}(\epsilon)\} < q_r$ and then $m$ so that the bound on $\delta_k$ and, hence, $\bar{\delta}$ is as small as needed.

5.3. *Generalizations.* 1. It seems likely that estimators which are adaptive with respect to $\sigma$ and $C$ could be constructed (for a fixed irrational $a$) by grouping frequencies $k$ within a given block $[q_n, q_{n+1})$ into a number of subblocks according to the value of $\|ka\|$ and then using some form of James–Stein shrinkage within each subblock. This methodology is now quite well established on other inverse problems with monotone eigenvalues; see, for example, Cavalier and Tsybakov (2002). Alternatively, adaptivity (up to logarithmic terms) is established via a wavelet deconvolution approach in Johnstone, Kerkyacharian, Picard and Raimondo (2004) for a class of Besov spaces including ellipsoids (13).

2. The ellipsoid results might also have been derived using the explicit evaluation of minimax risk given by Pinsker (1980). However, the method used here allows extension of the rate results to weighted $l_{2r}$ bodies of the form $\Theta = \{\theta : \sum k^{2\sigma r}\theta_k^{2r} \leq C^{2r}\}$ for $r \geq 1$ using essentially the same argument as for ellipsoids. For example, the analog of (58) states that if the ordered increasing $\epsilon_{(k)}$ corresponding to indices within a block $B_\nu$ satisfy some bound $\epsilon_{(k+1)}/\epsilon_{(k)} \leq \gamma$ (as happens for the boxcar $K_a$), then

$$b_\nu(t_\nu, \epsilon) = \sup\left\{\sum_{k \in B_\nu} \theta_k^2 \wedge \epsilon_k^2 : \sum_{k \in B_\nu} \theta_k^{2r} \leq t_\nu^{2r}\right\} \asymp \sum_1^{l_0} \epsilon_{(j)}^2,$$

where $l_0 = \sup\{l : \sum_{j=1}^l \epsilon_{(j)}^{2r} \leq t^{2r}\}$, and such sums can be estimated by the methods of this paper.

3. It is straightforward to extend the results of this paper to iterated kernels $K_a = ((2a)^{-1}I_{[-a,a]})^{\star m}$ with eigenvalues $r_k = (\sin \pi ka)^m/(\pi ka)^m$. However, kernels of the form $K_{a,b} = (2a)^{-1}I_{[-a,a]} \star (2b)^{-1}I_{[-b,b]}$ have eigenvalues

$$r_k = \frac{\sin \pi ka}{\pi ka}\frac{\sin \pi kb}{\pi kb} \asymp \frac{\|ka\|\|kb\|}{k^2 ab},$$

while the linear motion kernel (6) has

$$r_{k_1,k_2} = \frac{\sin \pi(k_1 a + k_2 ra)}{\pi(k_1 a + k_2 ra)}.$$

Considerable work exists on simultaneous Diophantine approximation problems [Schmidt (1980), Chapter 2], but whether this enables rate of convergence calculations is an open question.

## APPENDIX

PROOF OF (27) AND (28). We recall the convergence/divergence theorem of Khinchin [(1992), Theorem 32]. Let $\psi(x)$ be a positive continuous function of $x > 0$, such that $x\psi(x)$ is nonincreasing. Then the inequality $\|qa\| < \psi(q)$ has, for almost all $a$, a finite or infinite number of solutions in positive integers $q$ according as $\int_c^\infty \psi(x)\,dx$ converges or diverges.

For (27), consider $\psi(x) = (2x \log x)^{-1}$. Since the integral diverges, let $q$ be one of the infinitely many solutions to $\|qa\| < \psi(q)$ and choose $n$ so that $q_n \leq q < q_{n+1}$. It then follows from (22) and the property stated after (21) that

$$\frac{1}{2q_{n+1}} \leq \|q_n a\| \leq \|qa\| \leq \frac{1}{2q \log q} \leq \frac{1}{2q_n \log q_n},$$

from which (27) is immediate.

For (28), consider $\psi(x) = x^{-1}(\log x)^{-1-\delta}$. Since the integral converges, for all $q > q(a, \delta)$, we have $\|qa\| \geq \psi(q)$. In particular, from (22), for large $n$,

$$\frac{1}{q_{n+1}} \geq \|q_n a\| \geq \frac{1}{q_n (\log q_n)^{1+\delta}},$$

from which we obtain (28). □

PROOF OF (29). The method used to establish (29) for direct data may be extended in a straightforward manner to model (9), for example, by stepping through the arguments in Johnstone [(2003), Hyperrectangles chapter]. The key step in this approach, as in Donoho, Liu and MacGibbon (1990), is to establish that

$$R_L(\Theta, \epsilon) = \sup_{\tau \in \Theta} R_L(\Theta(\tau), \epsilon), \tag{67}$$

where $\Theta(\tau)$ is the hyperrectangle $\Pi[-\tau_i, \tau_i]$. This can be reduced to the Kneser–Kuhn minimax theorem [Johnstone (2003), Corollary A.4] applied to payoff function

$$f(c, s) = \sum_k [\epsilon^2 c_k^2 + (1 - c_k)^2 s_k], \tag{68}$$

defined for $(c, s) \in \ell_2(\mathbb{N}) \times \ell_1(\mathbb{N})$. But result (67) extends immediately to model (9) by replacing $\epsilon^2$ with $\epsilon_k^2$ in (68) and changing the domain of $c$ to the weighted Hilbert space $\ell_2(\mathbb{N}, (\epsilon_k^2)) = \{c : \sum c_k^2 \epsilon_k^2 < \infty\}$, and applying the minimax theorem in the same way. □

**Acknowledgments.** Iain M. Johnstone is grateful for hospitality of the Australian National University, where this project began. Marc Raimondo is grateful to Peter Hall for helpful references on Diophantine approximations. Part of this paper was written while Marc Raimondo visited Stanford University. Both authors are grateful for the comments of a referee, particularly for raising the issue discussed in Section 5.2.

## REFERENCES

BELITSER, E. and LEVIT, B. (1995). On minimax filtering over ellipsoids. *Math. Methods Statist.* **4** 259–273. MR1355248

BERTERO, M. and BOCCACCI, P. (1998). *Introduction to Inverse Problems in Imaging.* Institute of Physics, Bristol. MR1640759

CAVALIER, L., GOLUBEV, G. K., PICARD, D. and TSYBAKOV, A. B. (2002). Oracle inequalities for inverse problems. *Ann. Statist.* **30** 843–874. MR1922543

CAVALIER, L. and TSYBAKOV, A. B. (2002). Sharp adaptation for inverse problems with random noise. *Probab. Theory Related Fields* **123** 323–354. MR1918537

DONOHO, D. L. (1995). Nonlinear solution of linear inverse problems by wavelet-vaguelette decomposition. *Appl. Comput. Harmon. Anal.* **2** 101–126. MR1325535

DONOHO, D. L., LIU, R. C. and MACGIBBON, K. B. (1990). Minimax risk over hyperrectangles, and implications. *Ann. Statist.* **18** 1416–1437. MR1062717

EVANS, S. N. and STARK, P. B. (2002). Inverse problems as statistics. *Inverse Problems* **18** R55–R97. MR1929274

GOLUBEV, G. K. and KHAS'MINSKIĬ, R. Z. (1999). A statistical approach to some inverse problems for partial differential equations. *Problemy Peredachi Informatsii* **35** 51–66. MR1728907

GOLUBEV, G. and KHASMINSKII, R. (2001). A statistical approach to the Cauchy problem for the Laplace equation. In *State of the Art in Probability and Statistics*: *Festschrift for Willem R. van Zwet* 419–433. IMS, Beachwood, OH. MR1836573

GROENEBOOM, P. and JONGBLOED, G. (2003). Density estimation in the uniform deconvolution model. *Statist. Neerlandica* **57** 136–157. MR2035863

HALL, P., RUYMGAART, F., VAN GAANS, O. and VAN ROOIJ, A. (2001). Inverting noisy integral equations using wavelet expansions: A class of irregular convolutions. In *State of the Art in Probability and Statistics*: *Festschrift for Willem R. van Zwet* 533–546. IMS, Beachwood, OH. MR1836579

JOHNSTONE, I. M. (2003). Function estimation and Gaussian sequence models. Draft of a monograph. Available at www-stat.stanford.edu/˜imj.

JOHNSTONE, I. M., KERKYACHARIAN, G., PICARD, D. and RAIMONDO, M. (2004). Wavelet deconvolution in a periodic setting (with discussion). *J. R. Stat. Soc. Ser. B Stat. Methodol.* **66** 547–573, 627–652.

JOHNSTONE, I. M. and SILVERMAN, B. (1990). Speed of estimation in positron emission tomography and related inverse problems. *Ann. Statist.* **18** 251–280. MR1041393

KHINCHIN, A. Y. (1992). *Continued Fractions*. Dover, New York. MR1451873

KOO, J.-Y. (1993). Optimal rates of convergence for nonparametric statistical inverse problems. *Ann. Statist.* **21** 590–599. MR1232506

KOROSTELEV, A. and TSYBAKOV, A. (1993). *Minimax Theory of Image Reconstruction*. *Lecture Notes in Statist.* **82**. Springer, New York. MR1226450

KRESS, R. (1999). *Linear Integral Equations*, 2nd ed. Springer, New York. MR1723850

LANG, S. (1966). *Introduction to Diophantine Approximations*. Addison–Wesley, Reading, MA. MR209227

MAIR, B. and RUYMGAART, F. H. (1996). Statistical inverse estimation in Hilbert scales. *SIAM J. Appl. Math.* **56** 1424–1444. MR1409127

MATHÉ, P. and PEREVERZEV, S. V. (2001). Optimal discretization of inverse problems in Hilbert scales. Regularization and self-regularization of projection methods. *SIAM J. Numer. Anal.* **38** 1999–2021 (electronic). MR1856240

O'SULLIVAN, F. and ROY CHOUDHURY, K. (2001). An analysis of the role of positivity and mixture model constraints in Poisson deconvolution problems. *J. Comput. Graph. Statist.* **10** 673–696. MR1938974

PINSKER, M. (1980). Optimal filtering of square integrable signals in Gaussian white noise. *Problems Inform. Transmission* **16** 120–133. MR624591

SCHMIDT, W. (1980). *Diophantine Approximation*. *Lecture Notes in Math.* **785**. Springer, Berlin. MR568710

TENORIO, L. (2001). Statistical regularization of inverse problems. *SIAM Rev.* **43** 347–366 (electronic).
MR1861086

WAHBA, G. (1990). *Spline Models for Observational Data*. SIAM, Philadelphia. MR1045442

DEPARTMENT OF STATISTICS
STANFORD UNIVERSITY
STANFORD, CALIFORNIA 94305
USA
E-MAIL: imj@stanford.edu

SCHOOL OF MATHEMATICS
AND STATISTICS
CARSLAW BUILDING F07
UNIVERSITY OF SYDNEY
SYDNEY NSW 2006
AUSTRALIA
E-MAIL: marcr@maths.usyd.edu.au